\documentclass[runningheads]{llncs}
\usepackage{graphicx}
\usepackage{amsmath}
\usepackage{hyperref}
\usepackage{amssymb}
\usepackage{amsthm}
\usepackage{amsopn}
\usepackage{mathrsfs}
\usepackage{mathtools}
\usepackage{stmaryrd}
\usepackage[T1]{fontenc}
\usepackage{lipsum}
\usepackage{graphicx}
\graphicspath{ {/home/user/images/} } 
\usepackage{subcaption}

\usepackage{tikz}
  \usetikzlibrary{matrix}
  \usetikzlibrary{fit}
  \usetikzlibrary{backgrounds}
  \usetikzlibrary{arrows}
  \usetikzlibrary{shapes}
  \usetikzlibrary{decorations.pathmorphing}
\usepackage{tikz-cd}

\DeclareMathOperator{\End}{End}

\DeclareMathOperator{\GL}{GL}

\DeclareMathOperator{\id}{id}
\DeclareMathOperator{\vol}{vol}

\DeclareMathOperator{\Lie}{Lie}

\DeclareMathOperator{\U}{U}
\DeclareMathOperator{\Sp}{Sp}
\DeclareMathOperator{\Li}{Li}

\begin{document}
\newtheorem{defi}[theorem]{Definition}
\newtheorem{ex}[theorem]{Example}
\newtheorem{rem}[theorem]{Remark}
\newtheorem{lem}[theorem]{Lemma}
\newtheorem{prop}[theorem]{Proposition}

\title{
Gaussian distributions on Riemannian symmetric spaces in the large $N$ limit
\thanks{This work benefited from partial support by the European Research Council under the Advanced ERC Grant Agreement Switchlet n.670645. S.H. also received partial funding from the Cambridge Mathematics Placement (CMP) Programme. C.M. was supported by Fitzwilliam College and a Henslow Fellowship from the Cambridge Philosophical Society.}}

%
%

\author{Simon Heuveline\inst{1,2}
\and
Salem Said\inst{3}
\and
Cyrus Mostajeran\inst{4}}
%
\authorrunning{S. Heuveline et al.}
%
\institute{
Centre for Mathematical Sciences, University of Cambridge, United Kingdom
\email{sph48@cam.ac.uk}
\and 
Mathematical Institute, University of Heidelberg, Germany
\and
CNRS, University of Bordeaux, France
\and 
 Department of Engineering, University of Cambridge, United Kingdom
%
}
%
\maketitle              
\begin{abstract}
We consider Gaussian distributions on certain Riemannian symmetric spaces. In contrast to the Euclidean case, it is challenging to compute the normalization factors of such distributions, which we refer to as partition functions. In some cases, such as the space of Hermitian positive definite matrices or hyperbolic space, it is possible to compute them exactly using techniques from random matrix theory. However, in most cases which are important to applications, such as the space of symmetric positive definite (SPD) matrices or the Siegel domain, this is only possible numerically. Moreover, when we consider, for instance, high-dimensional SPD matrices, the known algorithms for computing partition functions can become exceedingly slow. Motivated by notions from theoretical physics, we will discuss how to approximate the partition functions in the large $N$ limit: an approximation that gets increasingly better as the dimension of the underlying symmetric space (more precisely, its rank) gets larger. We will give formulas for leading order terms in the case of SPD matrices and related spaces. Furthermore, we will characterize the large $N$ limit of the Siegel domain through a singular integral equation arising as a saddle-point equation.

\end{abstract}

\section{Gaussian distributions on Riemannian symmetric spaces}

It is widely known that an isotropic Gaussian distribution on $\mathbb{R}^N$ takes the form $p_{\mathbb{R}^N}(x; \bar{x}, \sigma) = (2\pi \sigma)^{-n/2}\exp(- \frac{1}{2\sigma^2}(x- \bar{x})^2)$, where $\bar{x}$ and $\sigma$ denote the mean and standard deviation, respectively. Gaussian distributions can naturally be generalized to Riemannian manifolds $(M,g)$ with the property
\begin{align}
    Z_M(\bar{x}, \sigma) \coloneqq \int_M \exp(- \frac{1}{2\sigma^2}d_g(x, \bar{x})^2)\;d\vol_g(x) <\infty,
\end{align}
where $d\vol_g(x)$ denotes the Riemannian volume measure on $M$ and $d_g: M\times M \rightarrow \mathbb{R}$ denotes the induced Riemannian distance function. We will refer to $Z_M$ as the \textbf{partition function}. In this case,
\begin{align}
   p_M(x;\bar{x}, \sigma)\coloneqq \frac{1}{Z_M(\bar{x}, \sigma)}\exp(-  \frac{1}{2\sigma^2}d_g(x, \bar{x})^2),
\end{align}
is a well-defined probability distribution on $M$.

In numerous applications requiring the statistical analysis of manifold-valued data, it is important to be able to compute the partition function $Z_M$. A key difference between $p_{\mathbb{R}^N}$ and $p_M$ in general is that $Z_M$ is typically a complicated integral while $Z_{\mathbb{R}^N}(\bar{x}, \sigma)= (2\pi \sigma)^{n/2}$ is readily available. On a general Riemannian manifold, $Z_M$ will depend in a highly non-linear and intractable way on the mean $\bar{x}$ as the dominant contribution to $Z_M$ will come from local data such as the curvature at $\bar{x}$. This complex dependence on local data is known from the heat-kernel approach to the Atiyah-Singer Index Theorem (see \cite{BGV13} chapters 2 and 4), for instance. Hence, to have a chance of actually computing the full partition function analytically, we should restrict to spaces that look the same locally at any point. Riemannian symmetric spaces formalize this intuition and we are fortunate that precisely such spaces appear in most applications of interest (see, for example, \cite{Ba2012,Fl2004,HSM21,Pennec2019,Sa16,Sa17,TSM21}).

\begin{defi}
A \textbf{Riemannian symmetric space} is a Riemannian manifold $(M,g)$ such that for each point $p \in M$ there is a global isometry $s_p:M \rightarrow M$ such that $s_p(p)=p$ and $d _p s_p=-\id_{T_pM}$, where $d_ps_p$ is the differential of $s_p$ at $p$ and $\id_{T_pM}$ denotes the identity operator on the tangent space $T_pM$.
\end{defi}
This definition emphasizes the geometric property that there are a lot of isometries on $M$. In fact, in the general theory of symmetric spaces \cite{Eb97,Hel79}, it can be shown that any symmetric space is isomorphic to a quotient of Lie groups $M\cong G/H$, where $H$ is a compact Lie subgroup of $G$, and that integrals over $M$ can be reduced to integrals over the respective Lie groups and their Lie algebras (more precisely, certain subspaces thereof) by the following proposition \cite{Sa17}.

\begin{prop}
\label{thm:int_sym}
Given an integrable function on a non-compact symmetric space $(M=G/H,g)$, we can integrate in the following way: 
\begin{align}
    \int_{M}f(x)d\vol_g(x)= C \int_H \int_{\mathfrak{a}}f(a,h)D(a)da dh
\end{align}
where $dh$ is the normalized Haar measure on $H$ and $da$ is the Lebesgue measure on a certain subspace $\mathfrak{a}\subset \mathfrak{g}=\Lie(G)$ \cite{Hel79,Sa17}. The function $D:\mathfrak{a}\rightarrow\mathbb{R}^+$ is given by the following product over roots $\lambda: \mathfrak{a}\rightarrow\mathbb{R}^+$ (with $m_{\lambda}$-dimensional root space) 
    $D(a)=\prod_{\lambda>0} \sinh^{m_\lambda}(|\lambda(a)|)$.
\end{prop}

Proposition \ref{thm:int_sym} and the fact that its isometry group acts transitively on a symmetric space allows us to prove that in fact, the partition function of a symmetric space does not depend on $\bar{x}$: $Z_{G/H}(\bar{x}, \sigma)= Z_{G/H}(\sigma)$ (see \cite{Sa17}).
Moreover, with Proposition \ref{thm:int_sym} at hand, the problem of determining $Z_{G/H}$ reduces to calculating an integral over $\mathfrak{a}$, which is a linear space. Computing the resulting integrals can now be achieved numerically using a specifically designed Monte-Carlo algorithm \cite{Sa16}.  Even though this works well for small dimensions, $\mathfrak{a}\cong \mathbb{R}^N$ is typically still a high-dimensional vector space when $G/H$ is high-dimensional. The known algorithms start to break down for $N\approx40 $ \cite{Sa17}, even though cases involving $N > 250$ are of relevance to applications such as electroencephalogram (EEG) based brain-computer interfaces \cite{Ba2012,Se17}. 
It should be noted as an aside, however, that for some spaces, $N$ does not depend on the dimension of the underlying symmetric space. For instance, in hyperbolic $d$-space, $N=1$ independently of $d$ \cite{HSM21}.

Fortunately, for the space of positive definite Hermitian matrices, the resulting integrals have previously been studied in the context of Chern-Simons theory \cite{Ma05,Ti04} and, in fact, they fit within the general theory of random matrices \cite{Me04}. In this paper, we will extend such approaches to study the corresponding integrals for a broader range of symmetric spaces of interest, including the space of symmetric positive definite (SPD) matrices and the Siegel domain.
We should like to note that we have recently discovered that ideas similar to the ones explored in this paper were developed independently by \cite{Ti20} from a slightly different perspective, which we will draw attention to throughout this work. However, we expect that our work will complement \cite{Ti20} by providing a stronger emphasis on the large $N$ limit, including an analysis of the saddle-point equation.

\section{Partition functions at finite $N$ and random matrices}

Motivated by applications, we will study the following spaces: 
    \begin{enumerate}
    \item
    \label{def:SPD}
    The spaces of symmetric, Hermitian and quaternionic Hermitian positive definite matrices will be denoted by $\mathcal{P}_{\mathbb{F}}(N)\cong \GL(N, \mathbb{F})/K$ where $K\in \{\U(N),\textrm{O}(N),\Sp(N)\}$ for $\mathbb{F}\in \{\mathbb{R}, \mathbb{C}, \mathbb{H}\}$, respectively. In each case, $\mathfrak{a} \subset \Lie(\GL(N, \mathbb{F}))\cong \End(N, \mathbb{F})$ is the $N$-dimensional space of diagonal matrices. We will denote the partition functions by $Z_{\beta}(\sigma)$ where $\beta\in \{1,2,4\}$, respectively. $\mathcal{P}(N)\coloneqq \mathcal{P}_{\mathbb{R}}(N)$ is commonly referred to as the \textbf{space of SPD matrices}.
    \item
    \label{def:Siegel}
    The \textbf{Siegel domain} $\mathbb{D}(N)\cong \Sp(2N, \mathbb{R})/\U(N)$ is the space of complex symmetric $N\times N$ matrices $z$ such that the imaginary part $\textrm{Im}(z)$ is a positive definite matrix. We will denote its partition function by $Z_S$.
\end{enumerate}

Using Proposition \ref{thm:int_sym}, it can be seen  \cite{HSM21} that their respective partition functions are given by

\begin{align}
\label{al:Z_beta}
    Z_{\beta}(\sigma)=\frac{C_{N,\beta}(\sigma)}{(2\pi)^NN!}  \int_{\mathbb{R}^N_+}\prod_{i=1}^N \Big(\exp\big(-\frac{\log^2(u_i)}{2\sigma^2}\big)\Big) |\Delta(u)|^{\beta} \prod_{i=1}^N du_i
\end{align}
where $\Delta(u)\coloneqq \prod_{i<j} (u_i-u_j)$ is the \textbf{Vandermonde determinant}, $C_{N,\beta}(\sigma)\coloneqq \frac{\omega_{\beta}(N) (2\pi)^N}{2^{NN_\beta}}\exp\Big(-NN_{\beta}^2 \frac{\sigma^2}{2}\Big)$ and $N_{\beta}\coloneqq \frac{\beta}{2}(N-1)+1$. This is in the well-known form of a random matrix partition function:
\begin{align}
\label{equ:generalpart}
Z(\sigma)=\int_{(a,b)^N}\prod_{i=1}^N \Big(\exp\big(-V(u_i;\sigma))\Big) |\Delta(u)|^{\beta} \prod_{i=1}^N du_i,
\end{align}
even though the potential takes the non-standard form $V_{SW}(x;\sigma)\coloneqq\frac{1}{2\sigma^2}\log^2(x)$ in our case. For instance, with the quadratic potential $V_Q(x;\sigma)\coloneqq \frac{1}{2\sigma^2}x^2$, the random matrix partition function from Equation \eqref{equ:generalpart} corresponds to the famous orthogonal, unitary and symplectic random matrix ensembles for $\beta \in \{1,2,4\}$, respectively \cite{Me04}, on which we will comment further in Example \ref{ex:Wig}. 
In fact, matrix models with the potential $V_{SW}$ also appear in the physics literature as the partition functions of $\U(N)$ Chern-Simons theory on $S^3$ (\cite{Ma05,Ti04}). 
This relation is not directly of use here, but it is worth mentioning that it provided our original inspiration for probing the structure behind the normalizing factors.
Moreover, the large $N$ limit of $\U(N)$ Chern-Simons theory is of physical interest and has been well studied in the theoretical and mathematical physics literature such as \cite{Ho03} Chapter 36.2 and \cite{AK07}, which motivates Section \ref{sec:large_N}.
$Z_S$ also turns out to be of the form \eqref{equ:generalpart} with $\beta=1$ and the potential $V_S(u;\sigma)\coloneqq \frac{\log^2(u+\sqrt{u^2-1})}{8\sigma^2}$:
\begin{align}
\label{eq:Z_S}
  Z_S(\sigma)=\vol (\U(N))2^{\frac{N(N+1)}{2}}N!\int_{(1,\infty)^N}
\prod \limits_{i=1}^N \exp \big(- V_{S}(u_i;\sigma)\big) \Delta(u) \prod \limits_{i=1}^N du_i.
\end{align}

Integrals of the form \eqref{equ:generalpart} can generally be calculated, by bringing the Vandermonde determinant to a suitable form involving orthogonal or skew-orthogonal polynomials.

\begin{defi}
\label{defi:ortho_pol}
Let $V:(a,b)\rightarrow \mathbb{R}$ be a given potential. A set of polynomials $\{R_i:i=1, \dots ,N\}$, with $R_j(x)=a_jx^j+ \dots$ being of degree $j$ is called 

\begin{enumerate}
    \item 
\textbf{orthogonal with potential $V$} if they form an orthonormal basis for the space of degree $N$ polynomials with respect to the inner product \begin{align}
   \langle f,g\rangle_2=\int_a^b \exp\big(-V(x)\big) f(x)g(x)dx.
\end{align}

\item
\textbf{skew-orthogonal with potential $V$} if they bring the skew-symmetric product  
\begin{align}
\label{eq:innerprd_1}
   \frac{1}{2} \langle f,g\rangle_1=\int_a^b \int_a^b f(x)g(y)\operatorname{sign}(x-y) \exp (-V(x)) \exp(-V(y))dxdy
\end{align}
to the standard form, meaning \begin{align}
\label{eq:standardform}
\begin{cases}
    \langle R_{2k},R_{2l} \rangle_1&=\langle R_{2k+1},R_{2l+1} \rangle_1=0\\
    \langle R_{2k},R_{2l+1} \rangle_1&=-\langle R_{2l+1},R_{2k} \rangle_1=\delta_{kl}.
\end{cases}
\end{align}
\end{enumerate}
\end{defi}

If we are given orthogonal or skew-orthogonal polynomials for a given potential $V$, then its corresponding partition function with $\beta=1$ or $\beta=2$, respectively, can be calculated in terms of the polynomials' leading order coefficients \cite{Me04}. A similar story, which we will not go into also works for the quaternionic case $\beta=4$.
Fortunately, orthogonal polynomials for the potential $V_{SW}$ have previously appeared in the literature as \textbf{Stieltjes-Wigert polynomials} \cite{Sz39}, which can be used to arrive at the following result \cite{HSM21}.

\begin{prop}
\label{thm:exactZ_2}
The partition function $Z_2$ for $\mathcal{P}_{\mathbb{C}}(N)$ is given by
\begin{align*}
Z_2(\sigma)= \frac{\omega_2(N)}{2^{N^2}}(2 \pi \sigma^2)^{\frac{N}{2}} \exp \Big((N^3-N)\frac{\sigma^2}{6}\Big)\prod \limits_{k=1}^{N-1}(1-e^{-k\sigma^2})^{N-k}.
\end{align*}
\end{prop}

\begin{rem}
It has come to our attention that very recently, there has been further progress in the case $\beta=2$. In the beautiful paper \cite{Fo20}, higher moments of the probability density are provided alongside further physical interpretations that go beyond the ones presented here.
Moreover, the eigenvalue density at finite $N$ is computed in \cite{Ti20}. 
\end{rem}

Unfortunately, even though \cite{Ti20} provides explicit calculations for the first few skew-orthogonal polynomials of $V_{SW}$, general expressions for such polynomials have yet to be found as noted in \cite{HSM21,Ti20}.
However, if we let $a_i$ and $b_i$ denote the leading order coefficients of skew-orthogonal polynomials $\{P_i=a_ix^i+\dots|i=1,\dots ,N\}$ for $V_{SW}$ and skew-orthogonal polynomials $\{Q_i=b_ix^i+ \dots|i=1,\dots ,N\}$ for $V_{S}$, respectively, then we can obtain exact formulas for $Z_1$ and $Z_S$ at finite $N$ in terms of $a_i$ and $b_i$ \cite{HSM21}. For technical convenience, we will focus on the even-dimensional cases ($N=2m$) in this work. The odd-dimensional cases can be treated in a similar manner subject to a number of technical modifications. See \cite{Ti20} for details on the treatment of the odd-dimensional case.

\begin{prop}
\label{thm:finiteN_SPD}
If $N=2m$, the partition functions $Z_1$ and $Z_S$ are given by 
 \begin{align}
    Z_1(\sigma)&=\frac{\omega_1(N)}{2^{Nm}} \exp\big(-N((N-1)/2+1)^2(\sigma^2/2)\big) \bigg (\prod_{l=1}^N a_l \bigg)^{-1} \\
    Z_S(\sigma)&=\vol (\U(N))2^{m^2(N+1)}\bigg (\prod_{l=1}^N b_l \bigg)^{-1}.
\end{align}
\end{prop}

The coefficients $a_i$ and $b_i$ can be found by bringing the inner product \eqref{eq:innerprd_1} into the standard form \eqref{eq:standardform}, which can be achieved numerically via a symplectic Gram-Schmidt algorithm \cite{Sa05}. Nonetheless, such a computation becomes exceedingly challenging for large $N$, which motivates the consideration of the large $N$ limit provided in the following section.

\section{Partition functions in the large $N$ limit and the saddle-point equation}
\label{sec:large_N}

We will now change our point of view and compute the large $N$ limit of the partition function $Z$, rather than trying to compute it at finite $N$. Here, the \textbf{large $N$ limit} 
specifically refers to the limit where $N \rightarrow \infty$ while the \textbf{`t Hooft parameter} $t\coloneqq \sigma^2 N$ is kept fixed.
In this limit, $Z$ has an asymptotic expansion in $N^{-2}$, known as the \textbf{genus expansion}  \cite{HSM21,Ma05}: 
\begin{equation}
\log(Z(\sigma))\sim \sum \limits_{g=0}^{\infty} f_g(t) N^{2-2g}.   
\end{equation} 
The first term in this series becomes an increasingly good approximation as $N$ gets larger, which is very useful since we are able to compute it.

Since fixing $t$ implies that $\sigma^2 \rightarrow 0$, the large $N$ limit is also referred to as the \textbf{double scaling limit} \cite{Ti20}. Moreover, in the physics literature, it is also known as the \textbf{planar limit} due to its interpretation in terms of planar Feynman diagrams \cite{HSM21} originally dating back to the 1970s where it has been studied in the context of quantum chromodynamics \cite{tH74}.
Interestingly, one may also consider other limits such as $\sigma^2 \rightarrow 0$ and $\sigma^2 \rightarrow \infty$ while $N$ is fixed, which we will not consider in this work (see \cite{Ti20} (II, 3. a, b) for further details).
In the case of $Z_2$, the large $N$ limit can be directly computed from Proposition \ref{thm:exactZ_2} and is found to be
\begin{align}
\label{eq:calculatian_Z_2}
    \frac{1}{N^2}\log (Z_2(\sigma))\sim - \frac{1}{2} \log \left(\frac{2N}{\pi}\right)+ \frac{3}{4}+ \frac{t}{6} - \frac{\Li_3(e^{-t})-\zeta(3)}{t^2}
\end{align}
where $\Li_3(x)\coloneqq \sum \limits_{k=1}^\infty \frac{x^k}{k^3}$ (for $|x|<1$) is the trilogarithm.

Equation \eqref{eq:calculatian_Z_2} crucially relies on having a closed form expression for the partition function $Z_2$ for any finite $N$ in the first place (Proposition \ref{thm:exactZ_2}). This is not the case for $Z_1$ and $Z_S$, since we do not know the asymptotics of $a_i$ and $b_i$ from Proposition \ref{thm:finiteN_SPD}. Therefore, we will now discuss a different, more powerful approach to the large $N$ limit inspired by ideas from theoretical physics: a saddle-point approximation.  This will enable us to directly obtain the large $N$ limit of $Z_2$ in Proposition \ref{thm:largeN} by solving a certain singular integral equation, the \textbf{saddle-point equation}.

As above, we split  $Z_{\beta}(\sigma)=C_{N,\beta}(\sigma)\Tilde{Z}_{\beta}(\sigma)$, where $\Tilde{Z}_{\beta}$ is in an appropriate form for the use of
saddle point methods as discussed in \cite{Ma05} (see Equation (1.46) therein):
\begin{align}
\label{eq:ZTilde}
 \Tilde{Z}_{\beta}(\sigma)=\frac{1}{N!}\int_{\mathbb{R}^N_+} \exp \big(N^2 \Tilde{V}_{SW}(\vec{\lambda},\beta;t) \big)\prod_{k=1}^N \frac{d\lambda_k}{2\pi},
\end{align}
where we have introduced the \textbf{effective potential}:
\begin{align}
\label{eq:S_eff}
    \Tilde{V}_{SW}\left(\vec{\lambda},\beta;t\right)\coloneqq-\frac{1}{2tN}\sum_{i=1}^N \log^2(\lambda_i)+ \frac{\beta}{N^2} \sum_{i<j} \log |\lambda_i-\lambda_j|.
\end{align}
In the large $N$ limit, the dominant contributions to $\Tilde{Z}_{\beta}$ will come from the saddle points of $\Tilde{V}_{SW}$. Note that this is analogous to a \textbf{semiclassical limit}, which is further discussed in \cite{HSM21}. For the Siegel domain, we have a similar effective potential, which takes the form
\begin{align}
\label{V_S_eff}
    \Tilde{V}_{S}\left(\vec{\lambda}; t\right) \coloneqq-\frac{1}{8tN}\sum_{i=1}^N \log^2\left(\lambda_i+\sqrt{\lambda_i^2-1}\right)+ \frac{1}{N^2} \sum_{i<j} \log |\lambda_i-\lambda_j|.
\end{align}

\begin{figure}
\centering

\includegraphics[width=\textwidth]{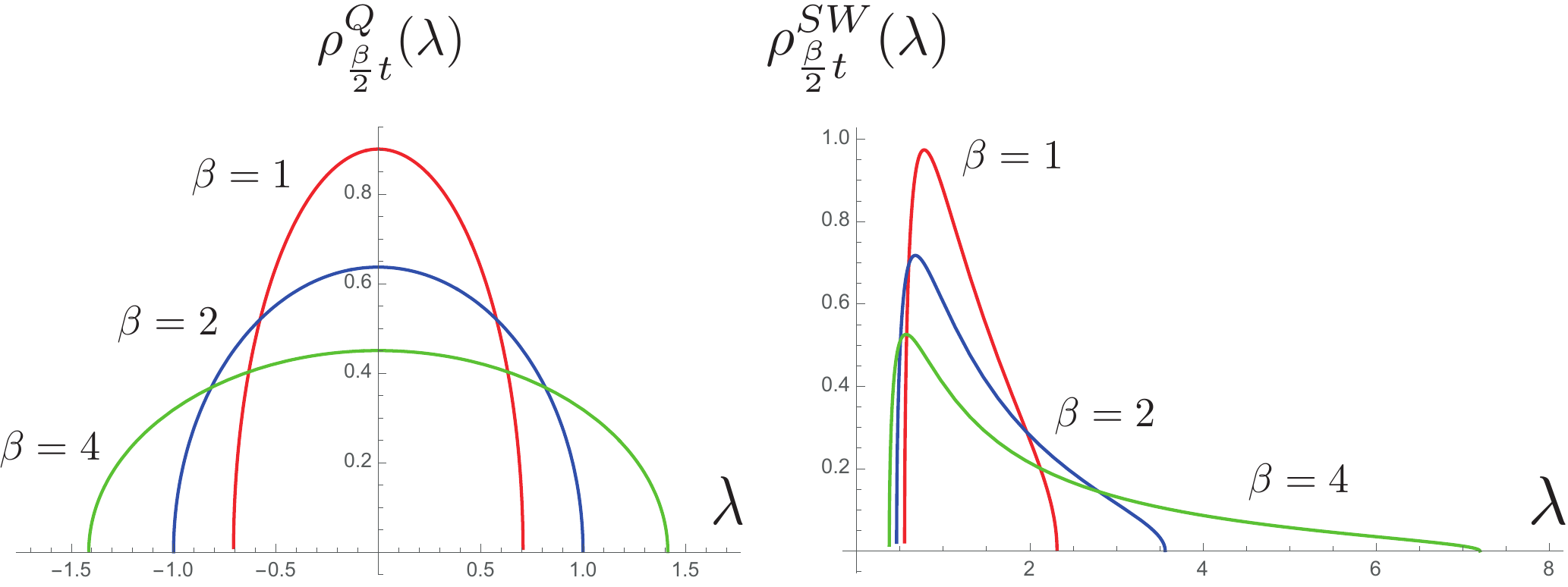}
    
    \caption{The master field $\rho^{Q}_{\beta t/2}$ of Wigner's semicircle law (left) and $\rho_{\beta t/2}^{SW}$ (right) for fixed $t=1/4$ and different choices of $\beta$. The red, blue and green curves correspond to the real, complex and quaternionic cases, respectively ($\beta = 1,2,4$).}

\label{fig:plots}

\end{figure}

\begin{rem}
\label{rem:semicirc}
 The effective potential $\Tilde{V}_{SW}$ gives rise to a physical interpretation of the theory and its large $N$ limit: the $N$ eigenvalues can be seen as static particles in the potential $V_{SW}$ interacting through a logarithmic Coulomb repulsion (the second term of the effective potential). As $\sigma^2=\frac{t}{N}$ decreases, the repulsion becomes weak and all particles can sit next to each other close to the minimum of the potential ($\lambda=1$), while the particles tend to spread out for large $\sigma^2$ as observed in Figure \ref{fig:plots}. Now, the large $N$ limit can be seen to correspond to the addition of more and more particles (i.e. $N\rightarrow \infty$) while letting their repulsion become increasingly weak (fixing $t=N\sigma^2$). Finding the limiting distribution $\rho_t$ (that still depends on $t$) turns out to characterize the large $N$ limit of the partition function. It can be obtained by solving the saddle-point equation for $\Tilde{V}$ in a continuum limit as discussed below. Motivated by the physics literature, we will refer to this $\rho_t$ as the \textbf{master field}. In the random matrix theory literature, this approach is known as the \textbf{Coulomb gas method} \cite{Me04}. 

\end{rem}

\begin{rem}
\label{rem:universality}
A further observation, from the effective action \eqref{eq:S_eff} is that up to an overall factor of $\Tilde{V}$ (which leaves the saddle-point equation invariant), we can absorb the $\beta$ by rescaling $t \mapsto t/\beta$. So, in the large $N$ limit, $\Tilde{Z}_{\beta}(\sigma) \sim \Tilde{Z}_{2}(\sqrt{\beta/2}\;\sigma)$ irrespective of the choice of $\beta \in \{1,2,4\}$, which is remarkable considering the quite distinct geometric origins associated with the different values of $\beta$.
 We refer to this phenomenon as \textbf{universality}: the three cases $\mathbb{K} \in \{\mathbb{R}, \mathbb{C}, \mathbb{H}\}$ start out differently, but "flow" to a universal limit (characterized by the master field) as $N\rightarrow \infty$. This is specially useful for the particularly important case of $\beta=1$, which can now be related to the large $N$ limit of $Z_2$ derived in Equation \eqref{eq:calculatian_Z_2}:
\begin{align}
    \frac{1}{N^2} \log (Z_{\beta}(\sigma))\sim \frac{1}{N^2}\left[ \log\left(Z_2\left( \sqrt{\tfrac{\beta}{2}} \sigma\right)\right) - \log\left(\frac{C_{\infty, \beta}(\sigma)}{C_{\infty,2}(\sqrt{\tfrac{\beta}{2}}\sigma)}\right)\right]
\end{align}
where $C_{\infty, \beta}$ are the large $N$ limits of the prefactors, which are discussed in \cite{HSM21}.
Below, we will give another formula for the large $N$ limit of $Z_\beta$, by solving the saddle-point equation explicitly.

\end{rem}

\begin{ex}[Wigner's semi-circle law]
\label{ex:Wig}
For simplicity, since $V_Q(\lambda;\sigma)=\frac{1}{2\sigma^2} \lambda^2$ has a simpler form than $V_{SW}$ or $V_S$ yet illustrates all the necessary ideas, we will begin by considering its saddle-point equation. As motivated in Remark \ref{rem:semicirc} and following  \cite{Ma05} ((1.47)-(1.53)), we are interested in the saddle points of the effective potential, which in this case reads 
    $\Tilde{V}_Q(\vec{\lambda};t)\coloneqq-\frac{1}{tN} \sum_{i} \lambda_i^2 + \frac{\beta}{N^2} \sum_{i<j} \log |\lambda_i-\lambda_j|$
in analogy with equations \eqref{eq:S_eff} and \eqref{V_S_eff}. The saddle points are characterized by $\frac{d}{d\lambda_k}\Tilde{V}_Q\left(\vec{\lambda};t\right)=0$ for all $k = 1, \dots, N$, which is simply \begin{align}
\label{eq:fin_saddle}
     \frac{1}{\beta t} \lambda_k=\frac{1}{N}\sum_{j\neq k} \frac{1}{\lambda_k-\lambda_j}=:P\left(\int \frac{\rho_{t,N}(\lambda')}{\lambda'-\lambda_k}d \lambda' \right) 
\end{align}
where $\rho_{t,N}$ is formally given by $\rho_{t,N}(\lambda)=\frac{1}{N} \sum_j \delta(\lambda-\lambda_j)$ and $P$ is a discrete Cauchy principal value.
 The large $N$ limit can now be regarded as a continuum limit, in which $\rho_{t,N}$ becomes a continuous function $\rho_t$. Equation \eqref{eq:fin_saddle} becomes 
\begin{align}
\label{eq:eom}
       \frac{1}{\beta t} \lambda= P\left(\int_{-\infty}^{\infty} \frac{\rho_t(\lambda')d \lambda'}{\lambda-\lambda'} \right)
       \iff      \frac{1}{2 t} \lambda= P\left(\int_{-\infty}^{\infty} \frac{\rho_{\frac{2}{\beta}t}(\lambda')d \lambda'}{\lambda-\lambda'} \right) 
\end{align}
where now $P$ is the actual Cauchy principal value. The saddle-point equation \eqref{eq:eom} can be solved using resolvent methods and in \cite{Ma05} (Equation (1.82)) it is 
shown that the solution in this case is simply a semicircle of radius $2 \sqrt{t}$: $ \rho^Q_t(\lambda)=\frac{1}{2\pi t} \sqrt{4t-\lambda^2} \chi_{\mathcal{C}(t)}(\lambda)$, where $\chi_{\mathcal{C}(t)}$ denotes the characteristic function supported on the interval 
$\mathcal{C}(t)=[-2\sqrt{t},2 \sqrt{t}]$. 
This celebrated semicircle law was first derived by Wigner in 1955 \cite{Wi55}.

\end{ex}

In fact, the steps of Example \ref{ex:Wig} that lead to the saddle point equation \eqref{eq:eom} can be performed analogously for the potentials $V_{SW}$ and $V_S$. Moreover, the saddle-point equation for $V_{SW}$ can even be solved analytically using resolvent methods \cite{HSM21,Ma05}, which leads to the following result.

\begin{prop}
\label{thm:largeN}
Let $\beta>0$. If $N\rightarrow\infty$ while the `t Hooft parameter $t=N\sigma^2$ is fixed, we get 
\begin{align}
    \frac{1}{N^2}\log( \Tilde{Z}_{\beta}(\sigma))\sim F_{uni}\left(\frac{\beta}{2} t\right)+ \mathcal{O}(N^{-2})
\end{align}
where \begin{align}
    F_{uni}(t)=-\frac{1}{2 t}\int_{\mathcal{C}(t)}\rho_t^{SW} (\lambda) \log^2\lambda \; d\lambda+ \int_{\mathcal{C}(t)^2} \rho_t^{SW}(\lambda) \rho_t^{SW} (\lambda') \log(|\lambda-\lambda'|) d\lambda d \lambda'
\end{align}
and the master field $\rho_t^{SW}$ is given by 
\begin{align}
\label{eq:rho}
   \rho_t^{SW} (\lambda)&= \frac{1}{\pi t \lambda} \tan^{-1}\left[ \frac{\sqrt{4 \lambda-(1+e^{-t}\lambda})^2}{1+e^{-t}\lambda} \right] \chi_{\mathcal{C}(t)}
\end{align}
where $\mathcal{C}(t)=[2e^{2t}-e^t+2e^{\frac{3t}{2}}\sqrt{e^t-1},2e^{2t}-e^t-2e^{\frac{3t}{2}} \sqrt{e^t-1}]$.
\end{prop}

In Figure \ref{fig:plots}, Wigner's semicircle distribution $\rho_t^Q$ is plotted alongside $\rho_t^{SW}$ for different choices of $\beta\in\{1,2,4\}$. It can be observed that as $t$ decreases (or equivalently, $\beta$ decreases), the distibutions tend to concentrate around the classical minima $\lambda=0$ (for $V_Q$) and $\lambda=1$ (for $V_{SW}$). Conversely, they tend to spread out as $t$ increases.

Finally, we note a similar equation that characterizes the master field in the case of the Siegel domain:
\begin{align}
\label{eq:saddle_S}
   \frac{\log\left(\lambda+\sqrt{\lambda^2-1}\right)}{4t\sqrt{\lambda^2-1}}=     P\left(\int_1^{\infty} \frac{\rho_t^S(\lambda')d \lambda'}{\lambda-\lambda'} \right).
\end{align}
At present, we are unaware of a solution to this equation and finding one either numerically or using resolvent methods is to be carried out in future work. Moreover, it is worth noting that we have certainly not exhausted the list of all (irreducible) symmetric spaces and have rather just focused on a few examples here. Even though further spaces are commented on in \cite{HSM21}, the development of a full description of the partition functions of all irreducible symmetric spaces and their interconnections remains to be done. Universalities and duality between compact and non-compact symmetric spaces will be guiding principles towards this goal.


\bibliographystyle{splncs04}
\bibliography{references}


\end{document}